\documentclass[reqno, 12pt]{amsart}
\usepackage{amssymb}
\usepackage{amsfonts}
\usepackage{enumerate}
\usepackage{cite}
\usepackage{amsmath}
\numberwithin{equation}{section}
\usepackage{amsfonts,amssymb,latexsym}
\usepackage[OT4]{fontenc}
\usepackage[utf8]{inputenc}
\usepackage{mathtools}
\usepackage{stmaryrd}
\usepackage{verbatim}
\usepackage{bbm}
\usepackage{dsfont}
\usepackage{amsthm}
\usepackage{etoolbox}
\usepackage[hyphens]{url}
\usepackage[hidelinks,hypertexnames=false]{hyperref}
\newcounter{localequation}

\AtBeginEnvironment{proof}{%
    \setcounter{localequation}{0}%
}

\newenvironment{localequation}
{%
    \refstepcounter{localequation}%
    \begin{equation*}
    \tag{\arabic{localequation}}
}
{%
    \end{equation*}
}

\newtheorem{theorem}{Theorem}[section]
\newtheorem{lemma}[theorem]{Lemma}
\newtheorem{proposition}[theorem]{Proposition}
\newtheorem{corollary}[theorem]{Corollary}
\theoremstyle{definition}
\newtheorem{definition}[theorem]{Definition}

\newtheorem*{conjecture}{Conjecture}
\theoremstyle{remark}
\newtheorem{remark}[theorem]{Remark}

\theoremstyle{definition}
\newtheorem{question}{Question}

\DeclareMathOperator{\diam}{diam}

\DeclareMathOperator{\Fix}{Fix}

\begin{document}
\title[Lau's Fixed Point Problem for Elementary Amenable Groups]{Lau's Fixed Point Problem for Elementary Amenable Groups and Finite Left Amenable Semigroups}
\author[S. Borzdy\'{n}ski]{S\l awomir Borzdy\'{n}ski}

\begin{abstract}
	Lau's fixed point conjecture is proved for elementary amenable groups.
	More generally, a weak$^*$-continuous nonexpansive representation of such a group on a nonempty weak$^*$-closed nonexpansive retract of a weak$^*$-compact convex subset $C$ of a dual Banach space has a nonempty common fixed point set that is a nonexpansive retract of $C$.
	For finite left amenable semigroups, continuity and weak$^*$-closedness of the retract are unnecessary.
	The argument also leads to a new class of amenable groups with the same retraction property, which contains all elementary amenable groups and may be strictly larger.
	This class may provide a useful framework for further research on Lau's conjecture.
	Bruck's theorem is also extended from abelian to elementary amenable group actions on nonexpansive retracts of weakly compact convex sets with the fixed point and conditional fixed point properties, without weak continuity, thus partially answering a longstanding question of Bruck.
	Several consequences of the theorem are recorded.
\end{abstract}

\address{S\l awomir Borzdy\'{n}ski, Institute of Computer Science and Mathematics, Maria Curie-Sk\l odowska University, 20-033 Lublin, Poland}
\email{slawomir.borzdynski@gmail.com}
\maketitle

\section{Introduction}

A fundamental result of Day asserts that every continuous affine action of a left amenable semigroup on a compact convex subset of a locally convex space has a common fixed point \cite{Day1961}.
Lau asked whether an analogous statement remains valid for jointly weak$^*$-continuous nonexpansive actions of left amenable semitopological semigroups on nonempty weak$^*$-compact convex subsets of dual Banach spaces \cite{Lau1976}.
This nonlinear fixed point conjecture remains open in its full generality.
Known positive cases include commuting families \cite{BorzynskiWisnicki2014}, norm-separable domains \cite[Theorem~6.2]{LauZhang2012}, and domains with the Radon--Nikod\'ym property \cite{Wisnicki2024}; see also \cite{LauZhang2018,LauZhangSurvey,LauZhang2022}.

The existence of a common fixed point is often only the first step.
A stronger conclusion is the existence of a nonexpansive retraction onto the common fixed point set.
The classical theorem of Bruck shows that this stronger conclusion can be obtained from an appropriate fixed point property \cite{Bruck1973}.
Further constructions for semigroup actions include the sunny nonexpansive retractions of Aleyner and Reich \cite{AleynerReich2005} and the ergodic retractions of Saeidi \cite{Saeidi2012}.
In fact Saeidi's retraction results play an essential role in the presented proofs multiple times.

The present work stems from the author's doctoral research~\cite{Borzynski2021Thesis}, in which it was shown that in the proof of the inductive retraction theorem from~\cite{Borzynski2018} commutativity can be replaced by a more general, albeit quite technical condition.
This observation motivated the search for an even broader and more natural criterion, which eventually led to the elementary amenable setting considered here.

In this paper, we prove Lau's conjecture for elementary amenable groups, regarded as discrete semitopological semigroups.
More precisely, we show that the common fixed point set of the representation is a nonempty nonexpansive retract of the ambient weak$^*$-compact convex set $C$, and for that, it is enough for the representation in question to be defined on a weak$^*$-closed nonexpansive retract, which is a subset of $C$.
No separability or countability assumptions are imposed.
As a by-product, the argument identifies a class of amenable groups that contains all elementary amenable groups.
This class is closed under quotients, extensions, and direct limits.
It is not known whether this class is strictly larger than the class of elementary amenable groups.

Using the same approach we extend Bruck's commuting-family theorem \cite{BruckCommonFixedPoint} to elementary amenable groups on nonexpansive retracts of weakly compact convex sets with the fixed point and conditional fixed point properties, without weak continuity.
In 1974, Bruck asked whether commutativity could be replaced by left reversibility \cite[Remark~3]{BruckCommonFixedPoint}.
Our result answers this longstanding question affirmatively for elementary amenable groups in the weakly compact case.
The author thanks Andrzej Wiśnicki of the University of Life Sciences in Lublin for suggesting that this approach might also be used to extend Bruck's theorem.

Using Saeidi's results \cite{Saeidi2012}, we obtain the ergodic retractions, which on convex domains with hereditary FPP can also be chosen to preserve closed convex invariant subsets.
A further consequence for reflexive Banach spaces follows from a preprint, not yet peer reviewed, made public while this paper was being finalized \cite{OpenAIReflexive2026}.

For finite left amenable semigroups, we obtain the analogous retraction result without continuity or weak$^*$-closedness assumptions on the retract.

We conclude the paper with questions and suggestions for further research.

\section{Preliminaries}

Let $S$ be a semigroup and let $(C,\tau)$ be a topological space.
A \emph{representation} or \emph{action} of $S$ on $C$ is a mapping
$$
	S \times C \to C, \qquad (s,x)\mapsto T_{s}x,
$$
that satisfies
$
	T_sT_tx=T_{st}x
$ for all $s,t\in S$ and $x\in C$.\footnote{If $S$ has an identity $e$, we do not require $T_e=\mathrm{id}_C$.}
In practice, we often identify a representation with the indexed family
$
	\{T_s\in C^C: s \in S\}
$
whose members satisfy the identity above.

A \emph{semitopological semigroup} $S$ is a semigroup equipped with a Hausdorff topology for which multiplication is separately continuous.
Recall Lau's conjecture; cf. \cite[Open Problem~2]{LauZhangSurvey}.
\begin{conjecture}[Lau]
	\label{conj:Lau}
	Every jointly weak$^*$-continuous representation of a left amenable semitopological semigroup by nonexpansive self-mappings of a nonempty weak$^*$-compact convex subset of a dual Banach space has a common fixed point.
\end{conjecture}
A discrete semigroup $S$ is \emph{left amenable} if there is a positive linear functional $\mu\colon\ell^\infty(S;\mathbb R)\to\mathbb R$ such that
$$
	\mu(\mathds{1})=1,\qquad \mu({}_{s}f)=\mu(f),\qquad ({}_{s}f)(t)=f(st)
$$
for every bounded real-valued $f$ and all $s,t\in S$.
For the semitopological version and background, see additionally \cite{DayAmenableSemigroups,PatersonAmenability,LauTakahashi1995}.
Unless stated otherwise, all groups and semigroups considered in this paper are equipped with the discrete topology.
This convention may seem to restrict our results to the discrete case of Lau's conjecture, but for elementary amenable groups it doesn't cause loss of generality.
Indeed, elementary amenability(recalled below) depends only on the underlying abstract group.
Hence our results will apply to semitopological groups whose underlying abstract groups are elementary amenable.

When $S$ is finite, a mean $\mu$ is a weighted average.
Indeed, put $p(a):=\mu(\mathds{1}_{\{a\}})$.
Since $f=\sum_{a\in S}f(a)\mathds{1}_{\{a\}}$, linearity, positivity and normalization give
$$
	\mu(f)=\sum_{a\in S}p(a)f(a),
	\qquad p(a)\geq0,
	\qquad \sum_{a\in S}p(a)=1.
$$
Thus $p$ has nonempty support.

For a family $\mathcal S$ of self-mappings of $C$, we write $\Fix\mathcal S$ for its common fixed point set, consisting of the points of $C$ fixed by every member of $\mathcal S$.
A subset $D\subseteq C$ is $\mathcal S$-invariant if $T_s(D)\subseteq D$ for every $s\in S$.
A nonempty subset $D\subseteq C$ is a nonexpansive retract of $C$ if there is a nonexpansive mapping $R\colon C\to D$ with $R|_D=\mathrm{id}_D$.

For a Banach space $E$, write $J\colon E\to E^{**}$ for the canonical isometric embedding, defined by $(Jx)(x^*)=x^*(x)$, and use $J^{-1}$ for its inverse on $J(E)$.
It identifies the weak topology on $E$ with the weak$^*$ topology on $J(E)$, thereby transferring compactness and convexity in both directions.
Since $J$ is a linear isometry, the correspondences $T\mapsto JTJ^{-1}$ and $T'\mapsto J^{-1}T'J$ likewise transfer nonexpansive mappings, representations and retractions back and forth between a subset $C\subseteq E$ and its image $J(C)$.
The common fixed point set of the conjugated representation is $J(\Fix\mathcal S)$.

The class $\mathrm{EG}$ of elementary amenable groups is the smallest class containing the finite and abelian groups and closed under subgroups, quotients, extensions, and direct limits \cite{Chou1980}.
We use the following Chou--Osin characterization \cite{Chou1980,Osin2002}; see also \cite[Proposition~3.3]{OzawaRordamSato2015}.
\begin{theorem}[Chou--Osin]
	\label{thm:chou-osin}
	Let $\mathrm{B}_0$ consist of the finite groups and the groups isomorphic to $\mathbb Z$.
	Then $\mathrm{EG}$ is the smallest class of groups containing the trivial group and closed under direct limits and extensions by groups from $\mathrm{B}_0$.
\end{theorem}
Here, extension by $\mathrm{B}_0$ means that $N\in\mathrm{EG}$ and $H\in\mathrm{B}_0$ imply $G\in\mathrm{EG}$ whenever $1\longrightarrow N\longrightarrow G\overset{\pi}{\longrightarrow}H\longrightarrow1$ is a short exact sequence.
In general for such sequences we identify $N$ with $\ker\pi\trianglelefteq G$, so $H\cong G/N$.
For the quotient map, $\pi(g)=\pi(h)$ precisely when $g=hn$ for some $n\in N$.
We also use $\operatorname{im}\varphi\cong G/\ker\varphi$ for a homomorphism $\varphi$.

A family of subgroups $(G_i)_{i\in I}$ is directed under inclusion if any two lie in some common $G_j$.
Its union $G:=\bigcup_{i\in I}G_i$ is then a subgroup, called the \emph{directed union} of the $G_i$.

A directed system $(G_i,\lambda_{ij})$ on a directed set $(I,\preceq)$ has homomorphisms $\lambda_{ij}\colon G_i\to G_j$ for $i\preceq j$, satisfying $\lambda_{ii}=\mathrm{id}_{G_i}$ and $\lambda_{jk}\lambda_{ij}=\lambda_{ik}$ for $i\preceq j\preceq k$.
Write its direct limit and canonical maps as
\begin{equation}
	\label{eq:direct-lim-hom}
	\lambda_i\colon G_i\to G=\varinjlim_jG_j,
	\qquad \lambda_j\lambda_{ij}=\lambda_i,\quad i\preceq j.
\end{equation}
We use the fact that every compatible family $\psi_i\colon G_i\to H$ factors uniquely as $\psi_i=\psi\lambda_i$ through a homomorphism $\psi\colon G\to H$.
Group theoretic facts stated here can be verified in \cite{DummitFoote2004,Rotman1995}.

\section{Results}
Throughout this section, $\tau$ is a Hausdorff vector topology on the ambient space.

\begin{definition}
\label{def:rf}
Let $C$ be a nonempty subset of a Banach space and let $\mathcal{A}$ be a class of semigroups.
For a Hausdorff topology $\tau$ on the ambient space, we say that $C$ has property $\mathcal{RF}_\tau(\mathcal{A})$ if, for every nonempty $\tau$-closed nonexpansive retract $D$ of $C$ and every representation $\mathcal S$ of $S\in\mathcal{A}$ by $\tau$-continuous nonexpansive self-mappings of $D$, the set $\Fix\mathcal S$ is a nonempty nonexpansive retract of $C$.
The retractions need not be $\tau$-continuous.

We write $\mathcal{RF}_0(\mathcal{A})$ for the same requirement with the $\tau$-closedness assumption on $D$ and the $\tau$-continuity assumption on the mappings omitted.
\end{definition}
\subsection{Lau's conjecture results}

First we adapt the approximation method of \cite[Theorem~3.1]{BorzynskiWisnicki2014} to mappings defined and continuous only on a retract.

\begin{lemma}[Single map]
	\label{lem:single-map-on-retract}
	Let $C$ be a nonempty norm-closed bounded convex $\tau$-compact subset of a Banach space whose norm is $\tau$-lower semicontinuous.
	Let $D$ be a nonempty $\tau$-closed nonexpansive retract of $C$ and let $T\colon D\to D$ be $\tau$-continuous and nonexpansive.
	Then $\Fix T$ is a nonempty nonexpansive retract of $C$.\footnote{
		Consequently, $C$ has property $\mathcal{RF}_\tau(\mathcal{A}_1)$, where $\mathcal{A}_1$ is the class of semigroups generated by one element.
	}
\end{lemma}

\begin{proof}
	Choose a nonexpansive retraction $R\colon C\to D$ and put $T':=TR$.
	Then $T'$ is nonexpansive, $T'(C)\subseteq D$, and $T'|_D=T$.
	For each $\alpha\in(0,1)$ and $x\in C$, the Banach contraction principle gives a unique point (cf. \cite{Bruck1973} or \cite[chap. 11]{GoebelKirk1990}) $B_\alpha x\in C$ satisfying
	\begin{localequation}
		\label{eq:single-map-implicit}
		B_\alpha x=(1-\alpha)x+\alpha T'B_\alpha x.
	\end{localequation}
	Subtracting the equations for $x,y\in C$ and using nonexpansiveness of $T'$ yields
	$$
		\|B_\alpha x-B_\alpha y\|
		\leq(1-\alpha)\|x-y\|+\alpha\|B_\alpha x-B_\alpha y\|.
	$$
	Thus $B_\alpha$ is nonexpansive, as is $P_\alpha:=T'B_\alpha\colon C\to D$.
	Since $T'|_D=T$, equation~\eqref{eq:single-map-implicit} and nonexpansiveness of $T'$ yield
	\begin{localequation}
		\label{eq:single-map-approx}
		\begin{aligned}
			\|TP_\alpha x-P_\alpha x\|
			&=\|{T'}^2B_\alpha x-T'B_\alpha x\|\leq\|T'B_\alpha x-B_\alpha x\|\\
			&=(1-\alpha)\|T'B_\alpha x-x\|\leq(1-\alpha)\diam C.
		\end{aligned}
	\end{localequation}

	Order $(0,1)$ as usual.
	Since $D$ is $\tau$-compact, Tychonoff's theorem gives compactness of $D^C$ in the topology of pointwise $\tau$-convergence.
	Thus a subnet $(P_{\alpha_\gamma})_\gamma$ converges pointwise in $\tau$ to a mapping $P\colon C\to D$.
	The $\tau$-lower semicontinuity of the norm gives
	$$
		\|Px-Py\|
		\leq\liminf_\gamma\|P_{\alpha_\gamma}x-P_{\alpha_\gamma}y\|
		\leq\|x-y\|,
		\qquad x,y\in C,
	$$
	so $P$ is nonexpansive.
	By the $\tau$-continuity of $T$ on $D$ and the approximate fixed-point estimate \eqref{eq:single-map-approx},
	$$
		\|TPx-Px\|
		\leq\liminf_\gamma\|TP_{\alpha_\gamma}x-P_{\alpha_\gamma}x\|
		=0.
	$$
	Hence $P(C)\subseteq\Fix T$.

	If $z\in\Fix T$, then $Tz=Rz=z$ and $T'z=z$.
	Uniqueness in~\eqref{eq:single-map-implicit} gives $B_\alpha z=z$, and hence $P_\alpha z=z$ for every $\alpha$.
	Passing to the limit gives $Pz=z$.
	Therefore $P$ is a nonexpansive retraction of $C$ onto $\Fix T$, which is nonempty because $C$ is nonempty.
\end{proof}

For commuting mappings, we pass the single map result to a finite family, and then use a compactness argument as in \cite[Theorem~2.9]{Borzynski2018} to obtain the general case.

\begin{lemma}[Commuting semigroups]
	\label{lem:commuting-family-on-retract}
	Let $C$ be as in Lemma~\ref{lem:single-map-on-retract}, and let $\mathcal{A}_{\mathrm c}$ be the class of commutative semigroups.
	Then $C$ has property $\mathcal{RF}_\tau(\mathcal{A}_{\mathrm c})$.
\end{lemma}

\begin{proof}
	Let $S\in\mathcal{A}_{\mathrm c}$, $D$ be a nonempty $\tau$-closed nonexpansive retract of $C$, and let $\mathcal S$ be a representation of $S$ by $\tau$-continuous nonexpansive self-mappings of $D$.
	For a finite subfamily $\mathcal U\subseteq\mathcal S$, the set $\Fix\,\mathcal U$ is $\tau$-closed in $C$, and commutativity gives $T(\Fix\,\mathcal U)\subseteq\Fix\,\mathcal U$ for every $T\in\mathcal S$.
	We show by induction on $|\mathcal U|$ that $\Fix\,\mathcal U$ is a nonempty nonexpansive retract of $C$.
	For $\mathcal U=\varnothing$, this is the assumption on $\Fix\varnothing=D$.
	If the assertion holds for $\mathcal U$ and we choose $T\in\mathcal S\setminus\mathcal U$, then $T|_{\Fix\,\mathcal U}$ is a $\tau$-continuous nonexpansive self-mapping of $\Fix\,\mathcal U$.
	Apply Lemma~\ref{lem:single-map-on-retract} to $T|_{\Fix\,\mathcal U}$ to obtain a nonexpansive retraction of $C$ onto
	$$
		\Fix(T|_{\Fix\mathcal U})=\Fix\,\mathcal U\cap\Fix T=\Fix(\mathcal U\cup\{T\}),
	$$
	completing the induction.

	Order the finite subfamilies of $\mathcal S$ by inclusion, and for each $\mathcal U$ choose a nonexpansive retraction $R_{\mathcal U}\colon C\to\Fix\,\mathcal U$.
	As in the proof of Lemma~\ref{lem:single-map-on-retract}, compactness of $D^C$ and the $\tau$-lower semicontinuity of the norm give a subnet $(R_{\mathcal U_\gamma})_\gamma$ converging pointwise in $\tau$ to a nonexpansive mapping $R\colon C\to D$.
	For a fixed $T\in\mathcal S$, eventually $T\in\mathcal U_\gamma$, so $R_{\mathcal U_\gamma}x\in\Fix T$ for every $x\in C$.
	Since $\Fix T$ is $\tau$-closed, it follows that $Rx\in\Fix T$.
	Hence $R(C)\subseteq\Fix\mathcal S$.
	Every $R_{\mathcal U}$ fixes $\Fix\mathcal S$, and therefore so does $R$.
	Thus $R$ is a nonexpansive retraction of $C$ onto $\Fix\mathcal S$, which is nonempty because $C$ is nonempty.
	Since $S$, $D$, and the representation $\mathcal S$ were arbitrary, $C$ has property $\mathcal{RF}_\tau(\mathcal{A}_{\mathrm c})$.
\end{proof}

The finite left amenable case combines DeMarr's argument \cite[Lemma~1]{DeMarr1963}, as used in \cite{HolmesLau1972}, with Saeidi's retraction theorem \cite{Saeidi2012}.
\begin{theorem}[Finite case]
	\label{lem:finite-lau}
	Let $\mathrm{LA}_{fin}$ be the class of finite left amenable semigroups.
	Every nonempty weak$^*$-compact convex subset $C$ of a dual Banach space has property $\mathcal{RF}_0(\mathrm{LA}_{fin})$.
	In particular, Lau's conjecture holds for every $S\in\mathrm{LA}_{fin}$.
\end{theorem}

\begin{proof}
	Let $S\in\mathrm{LA}_{fin}$ and let $\mathcal S=\{T_s:s\in S\}$ be a nonexpansive representation on a nonempty nonexpansive retract $D$ of $C$.
	Choose a nonexpansive retraction $R\colon C\to D$ and define
	$
	T'_s:=T_sR: C\to D
	$.
	Notice $RT'_t=T'_t$ and
	\begin{equation}
		\label{eq:retract-representation-composition}
		T'_sT'_t=
		T_sRT'_t=
		T_sT'_t=
		T_sT_tR=
		T_{st}R=
		T'_{st},
		\qquad s,t\in S.
	\end{equation}
	Thus $\mathcal S':=\{T'_s:s\in S\}$ is a nonexpansive representation on $C$.
	Moreover, $\Fix\mathcal S'=\Fix\mathcal S$, because every $T'_s$ takes values in $D$ and $T'_s|_D=T_s$.

	Let $\mu$ be a left invariant mean on $\ell^\infty(S;\mathbb R)$, with probability weights $p(a)=\mu(\mathds{1}_{\{a\}})$ as in the preliminaries.
	For $s,b\in S$, put $f=\mathds{1}_{\{b\}}$.
	Define ${}_{s}f(a)=f(sa)$.
	Using left invariance and Iverson notation\footnote{For a statement $P$, the Iverson bracket $\llbracket P\rrbracket$ equals $1$ if $P$ is true and $0$ otherwise.}, we have
	$$
		p(b)=\mu(f)=\mu({}_{s}f)=\sum_{a\in S}p(a)\llbracket sa=b\rrbracket.
	$$
	Consequently, the support $A:=\{a\in S:p(a)>0\}$ - which is nonempty - satisfies $sA=A$ for every $s\in S$.
	Indeed, positive $p$ at $a\in A$ implies positive $p$ at $sa\in sA$, while positive $p$ at $b\in A$ requires $a\in A$ with $sa=b$ and $p(a)>0$.

	Let $C_1\subseteq C$ be an arbitrary nonempty weak$^*$-closed convex $\mathcal S'$-invariant subset.
	By weak$^*$-compactness and the Kuratowski--Zorn lemma, $C_1$ contains a minimal nonempty weak$^*$-compact convex $\mathcal S'$-invariant subset $C_0$.
	Choose $y\in C_0$ and put
	$
		K:=\{T'_ay:a\in A\}.
	$
	The set $K$ is nonempty and finite, and
	$$
		T'_s(K)=\{T'_{sa}y:a\in A\}=K,
		\qquad s\in S.
	$$
	Suppose $d:=\operatorname{diam}K>0$.
	Since $K$ is finite and hence norm-compact, \cite[Lemma~2]{HolmesLau1972} gives $b\in\operatorname{co}K\subseteq C_0$ with $r:=\max_{y\in K}\|b-y\|<d$.
	As in the proof of \cite[Theorem~2]{HolmesLau1972}, put
	$
		K_0:=C_0\cap\bigcap_{y\in K}\overline B(y,r).
	$
	This set contains $b$ and is weak$^*$-compact and convex, since closed norm balls are weak$^*$-closed.
	Since $K$ is finite and $r<d$, some point of $K$ lies outside $K_0$, so $K_0\ne C_0$.
	For $z\in K_0$ and $s\in S$, fix $y\in K$ and choose $u\in K$ with $T'_su=y$, using $T'_s(K)=K$.
	Nonexpansiveness and $z\in K_0$ give $\|T'_sz-y\|=\|T'_sz-T'_su\|\leq\|z-u\|\leq r$.
	Since $y$ was arbitrary, $T'_sz$ lies in every $\overline B(y,r)$, $y\in K$, and also in $C_0$ by invariance of $C_0$.
	Thus $T'_sz\in K_0$, so $K_0$ is $\mathcal S'$-invariant, contradicting minimality of $C_0$.
	Therefore $K$ is a singleton.
	Its unique point is fixed by every $T'_s$ and belongs to $C_1$.
	Since $C_1$ was arbitrary, every nonempty weak$^*$-closed convex $\mathcal S'$-invariant subset of $C$ contains a common fixed point.
	By~\cite[Theorem~2.1(a)]{Saeidi2012}, there is a nonexpansive retraction of $C$ onto $\Fix\mathcal S'=\Fix\mathcal S$.
	Since $D$, $S$ and the representation were arbitrary, $C$ has property $\mathcal{RF}_0(\mathrm{LA}_{fin})$.
	Taking $D=C$ gives Lau's conjecture.
\end{proof}
Lau and Zhang \cite[Proposition~4.10]{LauZhang2012} obtain common fixed points for continuous actions on compact continuous retracts of their closed convex hulls, assuming the semigroup has the Schauder fixed point property.

Now we introduce an auxiliary class of groups.
Define $\mathcal{W}$ to be the largest class of discrete groups, with respect to inclusion, such that every nonempty weak$^*$-compact convex subset $C$ of any dual Banach space has property $\mathcal{RF}_{w^*}(\mathcal{W})$.

\begin{lemma}
\label{lem:extension-on-retract}
Let
$
1\longrightarrow N\longrightarrow G
\overset{\pi}{\longrightarrow} H\longrightarrow1
$
be a short exact sequence of groups.
If $N,H\in\mathcal{W}$, then $G\in\mathcal{W}$.
\end{lemma}

\begin{proof}
Let $C$ be a nonempty weak$^*$-compact convex subset of a dual Banach space $E$, and let $D$ be a nonempty weak$^*$-closed nonexpansive retract of $C$.
Consider a representation $\mathcal G=\{T_g:g\in G\}$ by weak$^*$-continuous nonexpansive self-mappings of $D$.
Since $N\in\mathcal{W}$, then definition of $\mathcal{W}$ implies that the set
$
C_0:=\Fix\{T_n:n\in N\}
$
is a nonempty nonexpansive retract of $C$.
Moreover, $C_0$ is weak$^*$-closed in $C$, since $D$ is weak$^*$-closed and every $T_n$ is continuous on $D$.

Normality of $N$ gives $T_g(C_0)\subseteq C_0$ for every $g\in G$.
Indeed, if $x\in C_0$ and $n\in N$, then
$$
T_nT_gx=T_gT_{g^{-1}ng}x=T_gx.
$$
Since $\pi$ is a surjective homomorphism, define $U_h:=T_g|_{C_0}$ for any $g\in\pi^{-1}(\{h\})$.
This is well-defined: if $\pi(a)=\pi(b)$, then $a=bn$ for some $n\in N$, so $T_ax=T_bT_nx=T_bx$ for every $x\in C_0$.
For $h,k\in H$, choose $a,b\in G$ with $\pi(a)=h$ and $\pi(b)=k$.
Then $\pi(ab)=hk$; we have $ab\in \pi^{-1}(hk)$ and
$$
U_hU_k=T_aT_b|_{C_0}=T_{ab}|_{C_0}=U_{hk}.
$$
Thus $\mathcal H:=\{U_h:h\in H\}$ is a representation of $H$ by weak$^*$-continuous nonexpansive self-mappings of $C_0$.
Since $H\in\mathcal{W}$, the set $\Fix\mathcal H$ is nonempty, nonexpansive retract of $C$.
And since $\Fix\mathcal H=\Fix\mathcal G$, we have $G\in\mathcal{W}$.
\end{proof}

\begin{lemma}
	\label{lem:grp-quotients}
	The class $\mathcal{W}$ is closed under surjective group homomorphisms.
	In particular, for $N\trianglelefteq G$
	we have
	$
		G\in\mathcal{W}\implies G/N\in\mathcal{W}.
	$
\end{lemma}

\begin{proof}
	Let $\pi\colon G\to H$ be a surjective homomorphism with $G\in\mathcal{W}$.
	Consider a representation $\{T_h:h\in H\}$ by weak$^*$-continuous nonexpansive self-mappings of a nonempty weak$^*$-closed nonexpansive retract $D$ of a nonempty weak$^*$-compact convex set $C$ in a dual Banach space.
	The family $\{T_{\pi(g)}:g\in G\}$ is a representation of $G$ on the same domain.
	Surjectivity gives
	$$
		\Fix\{T_{\pi(g)}:g\in G\}=\Fix\{T_h:h\in H\}.
	$$
	Since $G\in\mathcal{W}$, this common fixed point set is a nonempty nonexpansive retract of $C$.
	Thus $H\in\mathcal{W}$.
\end{proof}

\begin{lemma}
	\label{lem:direct-limit-canonical-images}
	Let $(G_i,\lambda_{ij})_{i\in I}$ be a directed system of groups, and let
	$
		G=\varinjlim_{i\in I}G_i
	$
	with the associated canonical homomorphisms\footnote{
		See equation~\eqref{eq:direct-lim-hom} in the preliminaries.
	}
	$
		\lambda_i\colon G_i\to G.
	$
	For every $i\in I$, put
	$
		L_i:=\lambda_i(G_i).
	$
	Then $(L_i)_{i\in I}$ is directed under inclusion and
	$$
		G=\bigcup_{i\in I}L_i.
	$$
\end{lemma}

\begin{proof}
	If $i\preceq j$, then $\lambda_i=\lambda_j\circ\lambda_{ij}$, so
	$L_i=\lambda_j(\lambda_{ij}(G_i))\subseteq L_j$.
	Thus $(L_i)$ is directed under inclusion, and its union $H:=\bigcup_iL_i$ is a subgroup of $G$.
	Define $\psi_i\colon G_i\to H$ by $\psi_i(g):=\lambda_i(g)$; this is well-defined since $\lambda_i(G_i)=L_i\subseteq H$.
	These homomorphisms satisfy $\psi_j\circ\lambda_{ij}=\psi_i$ whenever $i\preceq j$.
	The universal property gives a homomorphism $\psi\colon G\to H$ with $\psi\circ\lambda_i=\psi_i$ for every $i$.
	Let $\iota\colon H\to G$ be the inclusion.
	For every $i\in I$,
	$$
		(\iota\circ\psi)\circ\lambda_i
		=\iota\circ\psi_i
		=\lambda_i
		=\operatorname{id}_G\circ\lambda_i.
	$$
	Uniqueness in the universal property yields $\iota\circ\psi=\operatorname{id}_G$.
	Hence $\iota$ is surjective and $H=G$.
\end{proof}

Compare \cite[Lemma~2.1]{MonodTamuzTeryoshin2026} for the analogous closure under directed unions in the setting of continuous actions on compact convex spaces.

\begin{lemma}
	\label{lem:ur-directed-union}
	Let $G=\bigcup_{i\in I}G_i$ be a directed union of subgroups.
	If $G_i\in\mathcal{W}$ for every $i\in I$, then $G\in\mathcal{W}$.
\end{lemma}

\begin{proof}
	Let $C$ be a nonempty weak$^*$-compact convex subset of a dual Banach space, and let $D$ be a nonempty weak$^*$-closed nonexpansive retract of $C$.
	Consider a representation $\mathcal G=\{T_g:g\in G\}$ by weak$^*$-continuous nonexpansive self-mappings of $D$.
	Since $G_i\in\mathcal{W}$, choose a nonexpansive retraction $R_i\colon C\to K_i$, where $K_i:=\Fix\{T_g:g\in G_i\}$ is weak$^*$-closed in $D$.
	Order the net $(R_i)_i$ by inclusion of the subgroups $G_i$.
	As in the proof of Lemma~\ref{lem:commuting-family-on-retract}, a subnet $(R_{j_\gamma})_\gamma$ converges pointwise weak$^*$ to a nonexpansive mapping $R\colon C\to D$.
	For each fixed $i$, eventually $G_i\subseteq G_{j_\gamma}$, so $R_{j_\gamma}(C)\subseteq K_{j_\gamma}\subseteq K_i$.
	Closedness therefore gives $R(C)\subseteq\bigcap_i K_i=\Fix\mathcal G$.
	Since every $R_i$ fixes $\Fix\mathcal G$, so does $R$.
	Thus $R$ is a nonexpansive retraction of $C$ onto the nonempty set $\Fix\mathcal G$, and $G\in\mathcal{W}$.
\end{proof}

\begin{theorem}[Lau's conjecture for $\mathrm{EG}$]
	\label{thm:eg-on-retract}
	Let $C$ be a nonempty weak$^*$-compact convex subset of a dual Banach space $E$.
	Then $C$ has property $\mathcal{RF}_{w^*}(\mathrm{EG})$.
	In particular, Lau's conjecture holds for every elementary amenable group.
\end{theorem}

\begin{proof}
	Finite groups, including the trivial group, belong to $\mathcal{W}$ by Theorem~\ref{lem:finite-lau}, since $\mathcal{RF}_0$ implies $\mathcal{RF}_{w^*}$.
	Abelian groups, and in particular $\mathbb Z$, belong to $\mathcal{W}$ by Lemma~\ref{lem:commuting-family-on-retract}.
	Lemma~\ref{lem:extension-on-retract} shows that $\mathcal{W}$ is closed under extensions by groups from $\mathrm{B}_0$ (as defined in the preliminaries).

	To check closure under direct limits, let $G=\varinjlim G_i$ with $G_i\in\mathcal{W}$, and let $\lambda_i\colon G_i\to G$ be the canonical homomorphisms.
	By Lemma~\ref{lem:grp-quotients}, each image $L_i:=\lambda_i(G_i)$ belongs to $\mathcal{W}$.
	Lemma~\ref{lem:direct-limit-canonical-images} gives a directed union $G=\bigcup_i L_i$, so Lemma~\ref{lem:ur-directed-union} implies $G\in\mathcal{W}$.

	The Chou--Osin characterization, Theorem~\ref{thm:chou-osin}, now gives $\mathrm{EG}\subseteq\mathcal{W}$.
	By the definition of $\mathcal{W}$, $C$ has property $\mathcal{RF}_{w^*}(\mathcal{W})$, and hence $\mathcal{RF}_{w^*}(\mathrm{EG})$.
	Since $E$ and $C$ were arbitrary, taking $C$ itself as the retract domain proves Lau's conjecture for elementary amenable groups.
\end{proof}

\subsection{An extension of Bruck's theorem}

A nonempty bounded closed convex subset $C$ of a Banach space has the fixed point property (FPP) if every nonexpansive self-mapping has a fixed point, and the conditional fixed point property (CFPP) if each such mapping with a fixed point has one in every nonempty norm-closed convex invariant subset \cite{BruckCommonFixedPoint}.
Hereditary FPP requires FPP on every nonempty norm-closed convex subset of $C$ and implies both properties.

\begin{theorem}[Suggested by professor Wi\'snicki]
	\label{thm:bruck-eg-on-retract}
	Let $C$ be a nonempty weakly compact convex subset of a Banach space, and suppose that $C$ has FPP and CFPP.
	Then $C$ has property $\mathcal{RF}_0(\mathrm{EG})$.
\end{theorem}

\begin{proof}
	Let $\mathcal{M}$ be the largest class of groups such that every set $C$ satisfying the stated assumptions has property $\mathcal{RF}_0(\mathcal{M})$.

	Let $D$ be a nonempty nonexpansive retract of $C$, and consider a nonexpansive mapping $T\colon D\to D$.
	Choose a nonexpansive retraction $R\colon C\to D$ and put $T':=TR$.
	Since $C$ has FPP and CFPP and is weakly compact, \cite[Theorem~2]{BruckCommonFixedPoint} applied to $T'\colon C\to C$ yields a nonexpansive retraction of $C$ onto $\Fix T'$.
	Since $T'(C)\subseteq D$ and $R$ fixes $D$, $\Fix T'=\Fix T$.
	Thus $\Fix T$ is a nonempty nonexpansive retract of $C$.

	Now, for a representation $\mathcal{S}=\{T_n:n\in\mathbb Z\}$ on $D$, one has
	$
		\Fix\mathcal{S}=\Fix T_1.
	$
	Indeed, if $T_1x=x$, then $T_0x=T_0T_1x=T_1x=x$ and $T_{-1}x=T_{-1}T_1x=T_0x=x$, so every $T_n$ fixes $x$.
	Thus $\mathbb Z\in\mathcal{M}$.

	Finite groups, including the trivial group, are left amenable and belong to $\mathcal{M}$ by Theorem~\ref{lem:finite-lau}.
	Indeed, let $\mathcal S=\{T_g:g\in G\}$ be a nonexpansive representation of a finite group on $D$.
	The set $J(C)$ is weak$^*$-compact and convex.
	Recall that $J$ is an isometry, hence is injective and has an isometric inverse on $J(E)$.
	Thus $JRJ^{-1}\colon J(C)\to J(D)$ is nonexpansive and, for $x\in D$, $(JRJ^{-1})(Jx)=J(Rx)=Jx$.
	It is therefore a retraction onto $J(D)$.
	Apply Theorem~\ref{lem:finite-lau} to the conjugated representation $\{JT_gJ^{-1}\colon J(D)\to J(D)\}_{g\in G}$ to obtain a nonexpansive retraction $P'$ of $J(C)$ onto its common fixed point set.
	This set is $J(\Fix\mathcal S)$, since injectivity of $J$ gives, for $x\in D$ and $g\in G$,
	$$
		(JT_gJ^{-1})(Jx)=Jx
		\quad\Longleftrightarrow\quad J(T_gx)=Jx
		\quad\Longleftrightarrow\quad T_gx=x.
	$$
	Finally, $P:=J^{-1}P'J$ takes values in $\Fix\mathcal S$ and is nonexpansive since $J$ is an isometry.
	For $x\in\Fix\mathcal S$, we have $Px=J^{-1}P'(Jx)=J^{-1}Jx=x$, so $P$ is the required retraction.

	Now consider a short exact sequence
	$
		1\longrightarrow N\longrightarrow G\longrightarrow H\longrightarrow1
	$
	with $N,H\in\mathcal{M}$, and a representation $\mathcal G=\{T_g:g\in G\}$ on $D$.
	The property of $N$ yields a nonexpansive retraction of $C$ onto
	$
		C_0:=\Fix\{T_n:n\in N\}.
	$
	As in the proof of Lemma~\ref{lem:extension-on-retract}, normality gives $T_g(C_0)\subseteq C_0$, and the restrictions define a nonexpansive representation of $H$ on $C_0$ with common fixed point set $\Fix\mathcal G$.
	The property of $H$ supplies a retraction of $C$ onto that set.
	Thus $G\in\mathcal{M}$, and $\mathcal{M}$ is closed under extensions.

	The class $\mathcal{M}$ is closed under surjective homomorphic images.
	Indeed, let $\pi\colon G\to H$ be a surjective homomorphism with $G\in\mathcal{M}$, and let $\{T_h:h\in H\}$ be a nonexpansive representation on $D$.
	Since $\pi$ is a homomorphism, $\{T_{\pi(g)}:g\in G\}$ is a nonexpansive representation of $G$ on $D$.
	Since $\pi(G)=H$, this family consists of the same mappings as $\{T_h:h\in H\}$ and hence has the same common fixed point set.
	The property of $G$ gives a nonexpansive retraction of $C$ onto that set, so $H\in\mathcal{M}$.

	If $G=\bigcup_i G_i$ is a directed union with $G_i\in\mathcal{M}$, then the sets
	$
		D_i:=\Fix\{T_g:g\in G_i\}
	$
	are nonempty nonexpansive retracts of $C$ directed by reverse inclusion.
	By \cite[Lemma~3]{BruckCommonFixedPoint}, their intersection is a nonempty nonexpansive retract of $C$.
	This intersection equals $\Fix\mathcal G$, so $G\in\mathcal{M}$.
	Since Lemma~\ref{lem:direct-limit-canonical-images} reduces direct limits to surjective images and directed unions, the Chou--Osin characterization, Theorem~\ref{thm:chou-osin}, gives $\mathrm{EG}\subseteq\mathcal{M}$.
	Thus $C$ has property $\mathcal{RF}_0(\mathrm{EG})$.
\end{proof}

\begin{remark}
	This result extends the weakly compact case of \cite[Theorem~1]{BruckCommonFixedPoint} from abelian to elementary amenable groups: the common fixed point set is a nonempty nonexpansive retract of $C$.
	The domain of the representation may be any nonempty nonexpansive retract of $C$, without convexity or weak closedness; the geometric assumptions are required only for $C$.
	No weak continuity, separability or countability is required.
	This result partially answers \cite[Remark~3]{BruckCommonFixedPoint}.
\end{remark}

Combining Theorem~\ref{thm:bruck-eg-on-retract} with the invariant-mean construction of \cite[Lemmas~4.1--4.2]{Saeidi2012} gives invariant retractions.
On convex domains with hereditary FPP, \cite[Theorem~4.3(a)]{Saeidi2012} also gives preservation of norm-closed convex invariant subsets.

To facilitate the comparison of the next proof with Saeidi's method~\cite{Saeidi2012}, we conform to the notation used there: $\langle x^*,u\rangle:=u(x^*)$ for $x^*\in E^*$ and $u\in E^{**}$, and $\mu_g\varphi(g):=\mu(g\mapsto\varphi(g))$.
In particular, $\langle x^*,Jx\rangle=x^*(x)$.

\begin{corollary}
	\label{cor:bruck-eg-invariant-retraction}
	Let $C$ be as in Theorem~\ref{thm:bruck-eg-on-retract}, and let $\mathcal G=\{T_g:g\in G\}$ be a nonexpansive representation of $G\in\mathrm{EG}$ on a nonempty nonexpansive retract $D$ of $C$.
	Then there exists a nonexpansive retraction $P\colon C\to\Fix\mathcal G$ such that, for every $g\in G$,
	$$
		T_gP=P\quad\hbox{on }C,
		\qquad
		PT_g=P\quad\hbox{on }D.
	$$
	In particular, if $D=C$, then $T_gP=PT_g=P$ on $C$ for every $g\in G$.
\end{corollary}

\begin{proof}
	Choose nonexpansive retractions $R_D\colon C\to D$ and $R\colon C\to\Fix\mathcal G$, the latter supplied by Theorem~\ref{thm:bruck-eg-on-retract}.
	Put $T'_g:=T_gR_D$ for $g\in G$.
	As in equation~\eqref{eq:retract-representation-composition}, these mappings form a nonexpansive representation on $C$, with $\Fix\{T'_g:g\in G\}=\Fix\mathcal G$.

	Put $\mathcal{D}:=\ell^\infty(G;\mathbb R)$ and choose a right invariant mean $\mu\in \mathcal{D}^*$, which exists since $G$ is amenable.

	Consider the conjugated representation $\mathcal G_c:=\{JT'_gJ^{-1}:g\in G\}$ on $J(C)$, as in the proof of Theorem~\ref{thm:bruck-eg-on-retract}:
	$$
	    (JT'_gJ^{-1})(Jx)=J(T'_gx).
	$$
	Fix $x\in C$ and put $f_x(g):=J(T'_gx)$.
	Then $f_x(G)=\mathcal G_c(Jx)$; in other words, this is the orbit of $Jx$ under $\mathcal G_c$.
	Since $J$ is an isometry and $T'_gx\in C$, we have $\|f_x(g)\|=\|J(T'_gx)\|=\|T'_gx\|\leq\sup_{y\in C}\|y\|$; thus $f_x\colon G\to E^{**}$ is bounded.
	Moreover, the scalar function $g\mapsto\langle x^*,f_x(g)\rangle=x^*(T'_gx)$ belongs to $\mathcal{D}$ for every $x^*\in E^*$.
	Apply~\cite[Lemma~4.1]{Saeidi2012} with $E^{**}$ regarded as the dual space and $E^*$ as its predual.
	The lemma gives a unique element of $E^{**}$, which we denote by $F'_\mu(Jx)$, characterized by
	\begin{equation}
	    \label{eq:mean-value-definition}
	    \langle x^*,F'_\mu(Jx)\rangle
	    =\mu_g\langle x^*,f_x(g)\rangle,
	    \qquad x^*\in E^*.
	\end{equation}
	This element lies in the weak$^*$-closed convex hull of $f_x(G)$, and hence in $J(C)$, since $J(C)$ is weak$^*$-closed and convex.
	Since $x$ was arbitrary, this defines a mapping $F'_\mu\colon J(C)\to J(C)$.
	The pointwise inequalities $-\|\varphi\|_\infty\mathds{1}\leq\varphi\leq\|\varphi\|_\infty\mathds{1}$, together with positivity and $\mu(\mathds{1})=1$, give $|\mu(\varphi)|\leq\|\varphi\|_\infty$.
	For $x,y\in C$ and $x^*\in E^*$ with $\|x^*\|\leq1$, equation~\eqref{eq:mean-value-definition} and linearity of $\mu$ give
	$$
	    \begin{aligned}
	        |\langle x^*,F'_\mu(Jx)-F'_\mu(Jy)\rangle|
	        &=|\mu_g\langle x^*,f_x(g)-f_y(g)\rangle|\\
	        &\leq\sup_{g\in G}\|f_x(g)-f_y(g)\|
	        =\sup_{g\in G}\|T'_gx-T'_gy\|\\
	        &\leq\|x-y\|=\|Jx-Jy\|.
	    \end{aligned}
	$$
	Taking the supremum over these $x^*$ gives $\|F'_\mu(Jx)-F'_\mu(Jy)\|\leq\|Jx-Jy\|$, so $F'_\mu$ is nonexpansive.
	Set $F_\mu:=J^{-1}F'_\mu J\colon C\to C$, which is nonexpansive since $J$ is an isometry.
	For $x\in C$ and $x^*\in E^*$, equation~\eqref{eq:mean-value-definition} and the definition of $f_x$ give
	$$
	    \langle x^*,J(F_\mu x)\rangle
	    =\langle x^*,F'_\mu(Jx)\rangle
	    =\mu_g\langle x^*,f_x(g)\rangle
	    =\mu_g\langle x^*,J(T'_gx)\rangle.
	$$
	If $z\in\Fix\mathcal G$, then $f_z(g)=Jz$ for every $g$.
	By~\eqref{eq:mean-value-definition} and $\mu(\mathds{1})=1$, for every $x^*\in E^*$ we have
	$$
	    \langle x^*,F'_\mu(Jz)\rangle
	    =\mu_g\langle x^*,Jz\rangle
	    =\langle x^*,Jz\rangle\,\mu(\mathds{1})
	    =\langle x^*,Jz\rangle.
	$$
	This identity for all $x^*\in E^*$ gives $F'_\mu(Jz)=Jz$, and hence $F_\mu z=z$.
	Using $T'_gT'_h=T'_{gh}$ and right invariance of $\mu$, for $x\in C$ and $h\in G$ we obtain equality of the functionals $J(F_\mu T'_hx)$ and $J(F_\mu x)$ via the following identities for every $x^*\in E^*$:
	$$
	    \begin{aligned}
	        \langle x^*,J(F_\mu T'_hx)\rangle
	        &=\mu_g\langle x^*,J(T'_{gh}x)\rangle\\
	        &=\mu_g\langle x^*,J(T'_gx)\rangle
	        =\langle x^*,J(F_\mu x)\rangle.
	    \end{aligned}
	$$
	Injectivity of $J$ now yields $F_\mu T'_hx=F_\mu x$ for every $x\in C$, hence $F_\mu T'_h=F_\mu$ for every $h\in G$.

	Put $P:=RF_\mu$.
	This mapping is nonexpansive, takes values in $\Fix\mathcal G$, and fixes every point of $\Fix\mathcal G$, so it is a retraction onto that set.
	Its range gives $T_hP=P$ on $C$.
	For $x\in D$, we have $T'_hx=T_hx$ and therefore
	$$
	    PT_hx=RF_\mu T'_hx=RF_\mu x=Px.
	$$
	This proves the second identity on $D$.
\end{proof}

\begin{corollary}
	\label{cor:bruck-eg-preserving-retraction}
	Let $C$ be a nonempty weakly compact convex subset of a Banach space with hereditary FPP.
	For every representation $\mathcal G=\{T_g:g\in G\}$ of an elementary amenable group by nonexpansive self-mappings of $C$, there is a nonexpansive retraction $P\colon C\to\Fix\mathcal G$ such that, for every $g\in G$,
	$$
		T_gP=PT_g=P
	$$
	and $P(D)\subseteq D$ for every norm-closed convex $\mathcal G$-invariant subset $D$ of $C$.
\end{corollary}

\begin{proof}
	Conjugate the representation to $J(C)$ as in the proof of Theorem~\ref{thm:bruck-eg-on-retract}.
	Since $G$ is right amenable, we verify condition~(a) of \cite[Theorem~4.3]{Saeidi2012}.
	Let $D$ be an arbitrary nonempty norm-closed convex $\mathcal G$-invariant subset of $C$, as in the statement.
	Then $D$ is weakly compact, so $J(D)$ is weak$^*$-compact, hence weak$^*$-closed, convex and invariant in $J(C)$.
	Every nonempty weak$^*$-compact convex invariant subset of $J(C)$ arises in this way, since $J$ identifies the weak and weak$^*$ topologies.
	The set $D$ inherits hereditary FPP from $C$, and hence has both FPP and CFPP.
	Theorem~\ref{thm:bruck-eg-on-retract}, applied to the restricted representation on $D$ with retract domain $D$ itself, gives $D\cap\Fix\mathcal G\ne\varnothing$.
	Thus every such subset of $J(C)$ contains a common fixed point of the conjugated representation, including $J(C)$ itself.

	All the assumptions of \cite[Theorem~4.3(a)]{Saeidi2012} are now satisfied, so it gives a nonexpansive retraction $P'\colon J(C)\to J(\Fix\mathcal G)$.
	Then
	$$
		P:=J^{-1}P'J\colon C\to\Fix\mathcal G
	$$
	is a nonexpansive retraction with $T_gP=PT_g=P$ for all $g\in G$.
	The preservation assertion of the same theorem gives $P'(J(D))\subseteq J(D)$, and hence $P(D)\subseteq D$, for every $D$ in the statement.
\end{proof}

The next corollary follows from OpenAI's recently announced solution of the longstanding fixed point problem for reflexive Banach spaces \cite[Theorem~1.1]{OpenAIReflexive2026}.
At the time of writing, this preprint has not undergone peer review and the way it was presented to the community alongside other announced solutions sparked controversy.

\begin{corollary}
	\label{cor:bruck-eg-reflexive}
	Let $E$ be a reflexive Banach space and let $C$ be a nonempty norm-closed bounded convex subset of $E$.
	Then $C$ has property $\mathcal{RF}_0(\mathrm{EG})$.
	In particular, for every nonempty nonexpansive retract $D$ of $C$ and every representation $\mathcal G=\{T_g:g\in G\}$ of an elementary amenable group by nonexpansive self-mappings of $D$, the set $\Fix\mathcal G$ is nonempty and there is a nonexpansive retraction $P\colon C\to\Fix\mathcal G$ such that, for every $g\in G$,
	$$
		T_gP=P\quad\hbox{on }C,
		\qquad
		PT_g=P\quad\hbox{on }D.
	$$
	If $D=C$, then $P$ can additionally be chosen to preserve every norm-closed convex $\mathcal G$-invariant subset of $C$.
\end{corollary}

\begin{proof}
	By \cite[Theorem~1.1]{OpenAIReflexive2026}, every nonexpansive self-mapping of every nonempty norm-closed convex subset $D$ of $C$ has a fixed point, since $D$ is bounded and $E$ is reflexive.
	Thus $C$ has hereditary FPP, and hence both FPP and CFPP.
	Reflexivity also implies that $C$ is weakly compact.
	Theorem~\ref{thm:bruck-eg-on-retract} gives $\mathcal{RF}_0(\mathrm{EG})$ for $C$ and the nonempty common fixed point set.
	Corollary~\ref{cor:bruck-eg-invariant-retraction} gives the retraction with the stated identities.
	When $D=C$, Corollary~\ref{cor:bruck-eg-preserving-retraction} supplies a retraction with the same identities that also preserves every norm-closed convex $\mathcal G$-invariant subset of $C$.
\end{proof}

\subsection{Further directions}

Every group in $\mathcal{W}$ is amenable, since its fixed point property implies left amenability \cite[p.~530]{LauTakahashi2001}.
Motivated by this observation, the author proposes calling the groups in $\mathcal{W}$ \emph{retractively amenable} and denoting this class by $\mathrm{RA}$.

The class $\mathrm{RA}$ is closed under surjective homomorphic images (Lemma~\ref{lem:grp-quotients}), extensions (Lemma~\ref{lem:extension-on-retract}), and direct limits (proof of Theorem~\ref{thm:eg-on-retract}).
Together with the amenability observation above, the last paragraph of the proof of Theorem~\ref{thm:eg-on-retract} gives
$\mathrm{EG}\subseteq\mathrm{RA}\subseteq\mathrm{AG}$,
where $\mathrm{AG}$ denotes the class of amenable groups.

We also record the following closure property.
\begin{proposition}
\label{lem:ra-subgroups}
The class $\mathrm{RA}$ is closed under subgroups.
Consequently, a group belongs to $\mathrm{RA}$ if and only if all its finitely generated subgroups belong to $\mathrm{RA}$.
\end{proposition}
\begin{proof}
Let $H\leq G\in\mathrm{RA}$, and consider arbitrary $C\subseteq E^*$, a retract $D$ and a representation $\mathcal H=\{T_h:h\in H\}$ on $D$ as in the definition of $\mathrm{RA}$ and choose a nonexpansive retraction $R\colon C\to D$. We will construct an appropriate representation of $G$.

To this end set $C':=C^G$.
For $f\in C'$ we will write $f_t:=f(t)$.
Of course by Tychonoff's theorem $C'$ is compact in the product weak$^*$ topology.
But in order to use property $\mathrm{RA}$ we need $C'$ to be weak$^*$ compact in the topology of some dual space.

Notice that since $C$ is norm-bounded, $C'$ is a bounded convex subset of $\ell^\infty(G,E^*)=\ell^1(G,E)^*$.
Choose a closed norm ball $B$ in this dual space which is centered at $0$ and contains $C'$.
Let $\tau_1$ be the weak$^*$ topology on $B$ induced by $\ell^1(G,E)$, and let $\tau_2$ be the topology on $B$ inherited from the product weak$^*$ topology.
By the Banach--Alaoglu theorem, $(B,\tau_1)$ is compact, while $(B,\tau_2)$ is Hausdorff since the weak$^*$ topology on $E^*$ is Hausdorff.

By the definition of $\tau_1$, for every fixed $u\in\ell^1(G,E)$ the scalar map
$
    f\mapsto\langle f,u\rangle
$
is continuous from $(B,\tau_1)$ to $\mathbb R$.
Fix $t\in G$ and consider the coordinate projection $\pi_t\colon(B,\tau_1)\to(E^*,w^*)$, $\pi_t(f)=f_t$.
For each $x\in E$, choose $u\in\ell^1(G,E)$ given by $u_s=\llbracket s=t\rrbracket x$ for $s\in G$, using Iverson brackets.
Then
$
    \langle f,u\rangle=f_t(x)=(\pi_t(f))(x),
$\footnote{
    Under the identification $\ell^\infty(G,E^*)=\ell^1(G,E)^*$, a bounded family $f=(f_s)_{s\in G}$ represents the bounded linear functional
    $u\mapsto\langle f,u\rangle:=\sum_{s\in G}f_s(u_s)$.
}
so the scalar map $f\mapsto(\pi_t(f))(x)$ is continuous on $(B,\tau_1)$.
Since this holds for every $x\in E$, the definition of the weak$^*$ topology on $E^*$ implies that $\pi_t$ is continuous.

Now the identity map
$
    \operatorname{id}\colon(B,\tau_1)\longrightarrow(B,\tau_2)
$
is also continuous, since $\pi_t\circ\operatorname{id}=\pi_t$ for every $t\in G$.
It is a bijection from a compact space to a Hausdorff space, hence a homeomorphism.
The continuity of the identity and its inverse gives $\tau_1=\tau_2$, so their restrictions to $C'$ agree as well.
Hence $C'$ is weak$^*$-compact in $\ell^1(G,E)^*$.

Now, as in the proof of \cite[Proposition~2.2]{MonodTamuzTeryoshin2026}, put
\begin{localequation}
    \label{eq:ra-subgroups-equivariant}
    D':=\{f\in D^G:f_{ht}=T_hf_t\text{ for all }h\in H,\ t\in G\}.
\end{localequation}
This set is weak$^*$-closed, since $D$ is weak$^*$-closed and the mappings $T_h$ are weak$^*$-continuous.
Choose a set $A\subseteq G$ containing exactly one representative of each right coset of $H$ in $G$.
Every $t\in G$ then has a unique expression $t=\eta(t)\sigma(t)$, which defines maps $\eta\colon G\to H$ and $\sigma\colon G\to A$.
Uniqueness gives $\sigma(ht)=\sigma(t)$ and $\eta(ht)=h\eta(t)$ for $h\in H$ and $t\in G$.
Define $R'\colon C'\to D'$ by
$$
    (R'f)_t:=T_{\eta(t)}R(f_{\sigma(t)}).
$$
Using these identities and the representation property, we obtain
$$
    \begin{aligned}
        (R'f)_{ht}
        &=T_{\eta(ht)}R(f_{\sigma(ht)})=T_{h\eta(t)}R(f_{\sigma(t)})\\
        &=T_hT_{\eta(t)}R(f_{\sigma(t)})=T_h(R'f)_t,
    \end{aligned}
$$
so $R'$ is well-defined as a map into $D'$.
Nonexpansiveness of $R$ and the $T_h$ gives $\|R'f_1-R'f_2\|_\infty\leq\|f_1-f_2\|_\infty$ for $f_1,f_2\in C'$.
For $f\in D'$, we have
$$
    (R'f)_t=T_{\eta(t)}R(f_{\sigma(t)})=T_{\eta(t)}f_{\sigma(t)}=f_{\eta(t)\sigma(t)}=f_t.
$$
Thus $D'$ is a nonempty nonexpansive retract of $C'$.

For $g\in G$, define $T'_g\colon D'\to D'$ by
\begin{localequation}
    \label{eq:ra-subgroups-shift}
    (T'_gf)_t:=f_{tg}.
\end{localequation}
By~\eqref{eq:ra-subgroups-equivariant}, these shifts indeed map $D'$ into itself, since
$$
    (T'_gf)_{ht}=f_{htg}=T_hf_{tg}=T_h(T'_gf)_t.
$$
Moreover,
$$
    (T'_gT'_kf)_t=(T'_kf)_{tg}=f_{tgk}=(T'_{gk}f)_t.
$$
For fixed $g$, the map $t\mapsto tg$ is a bijection of $G$.
Consequently, for $f_1,f_2\in D'$,
$$
    \begin{aligned}
        \|T'_gf_1-T'_gf_2\|_\infty
        &=\sup_{t\in G}\|(f_1)_{tg}-(f_2)_{tg}\|\\
        &=\|f_1-f_2\|_\infty.
    \end{aligned}
$$
For each $t\in G$, the map $\pi_t\circ T'_g=\pi_{tg}|_{D'}$ is continuous as the restriction of a coordinate projection.
As above, continuity of all coordinate maps implies that $T'_g$ is continuous in the product weak$^*$ topology.
By the equality of topologies on $C'$ established above, it is also continuous in the weak$^*$ topology induced by $\ell^1(G,E)$.
Thus $\mathcal G':=\{T'_g:g\in G\}$ is a representation by weak$^*$-continuous isometries of $D'$.
Since $G\in\mathrm{RA}$, and all prerequisites on $C'$, $D'$ and $\mathcal G'$ are met, there is a nonexpansive retraction $P'\colon C'\to\Fix\mathcal G'$.

Define
\begin{localequation}
    \label{eq:ra-subgroups-diagonal-evaluation}
    \begin{aligned}
        \Delta&\colon C\to C', & (\Delta x)_t&:=x,\\
        \nabla&\colon C'\to C, & \nabla(f)&:=f_{e_G}.
    \end{aligned}
\end{localequation}
Thus $\Delta x\colon G\to C$ is the constant function with value $x$.
By~\eqref{eq:ra-subgroups-diagonal-evaluation}, the map $\Delta$ is isometric and $\nabla$ is nonexpansive.
If $f\in\Fix\mathcal G'$, then~\eqref{eq:ra-subgroups-shift} gives
$$
    f_g=(T'_gf)_{e_G}=f_{e_G},\qquad g\in G,
$$
so $f=\Delta x$ with $x:=f_{e_G}\in D$.
Then taking $t=e_G$ in~\eqref{eq:ra-subgroups-equivariant} gives $T_hx=T_hf_{e_G}=f_h$.
Since $f$ is constant with value $x$, we obtain $T_hx=x$ for every $h\in H$, hence $x\in\Fix\mathcal H$.
Conversely, if $x\in\Fix\mathcal H$, the constant function $\Delta x$ belongs to $D'$ and is fixed by every shift $T'_{g}$.
Therefore $\Fix\mathcal G'=\Delta(\Fix\mathcal H)$.
It follows that
$$
    P:=\nabla P'\Delta\colon C\to\Fix\mathcal H,
    \qquad Px=\bigl(P'(\Delta x)\bigr)_{e_G},
$$
is nonexpansive.
For $x\in\Fix\mathcal H$, we have $\Delta x\in\Fix\mathcal G'$, so $P'(\Delta x)=\Delta x$ and $Px=\nabla\Delta x=x$.
Thus $P$ is the required retraction, proving $H\in\mathrm{RA}$.

The forward implication of the final assertion follows from the subgroup closure just proved.
Conversely, suppose all finitely generated subgroups of $G$ belong to $\mathrm{RA}$.
They form a directed family.
Their union is $G$, since every $g\in G$ belongs to $\langle g\rangle$.
Lemma~\ref{lem:ur-directed-union} therefore yields $G\in\mathrm{RA}$.
\end{proof}

Recall that the class $\mathrm{SA}$ of subexponentially amenable groups satisfies
$
	\mathrm{EG}\subsetneq \mathrm{SA}\subsetneq \mathrm{AG}.
$
The first strict inclusion is witnessed by the first Grigorchuk group \cite{Chou1980,Grigorchuk1985}.
The Basilica group is not subexponentially amenable \cite[Proposition~13]{GrigorchukZuk2002} but is amenable \cite{BartholdiVirag2005}, and hence witnesses the second strict inclusion.

\begin{question}
	How is $\mathrm{RA}$ positioned relative to the classes of elementary amenable and subexponentially amenable groups?
	More specifically, is the inclusion $\mathrm{EG}\subseteq\mathrm{RA}$ strict, and does either of the inclusions
	$
		\mathrm{SA}\subseteq \mathrm{RA}
	$
	and
	$
		\mathrm{RA}\subseteq\mathrm{SA}
	$
	hold?
\end{question}

The author knows no example in $\mathrm{RA}\setminus\mathrm{EG}$ and does not know whether the first Grigorchuk group belongs to $\mathrm{RA}$.
Since every group is the directed union of its finitely generated subgroups, closure of $\mathrm{EG}$ under directed unions and subgroup closure of $\mathrm{RA}$ show that $\mathrm{RA}\setminus\mathrm{EG}$, if nonempty, contains a finitely generated group.

The definition of retractive amenability extends directly to discrete semigroups: call a semigroup $S$ \emph{retractively amenable} if every nonempty weak$^*$-compact convex subset $C$ of any dual Banach space has property $\mathcal{RF}_{w^*}(\{S\})$.
The resulting class consists of left amenable semigroups by the aforementioned fixed point characterization, and its restriction to groups recovers $\mathrm{RA}$.
This fact and Lau's original semigroup setting suggest the following research direction.

\begin{question}
	Can one define a natural, conservative\footnote{
        Here, \emph{conservative} is meant in the sense that restricting such definitions to groups recovers the corresponding class $\mathrm{EG}$.
	}
	semigroup analogue of elementary amenable groups?
	Would the resulting classes of semigroups be left amenable?
	Would the analog of the relation $\mathrm{EG}\subseteq\mathrm{RA}$ still follow?
\end{question}
A rough idea here is to adapt the operations in the Chou--Osin characterization to semigroups.

Notice that the definition of retractive amenability is analytic in nature:
it is formulated in terms of nonexpansive
retractions.
This leads to yet another question.

\begin{question}
	Is there an intrinsic group-theoretic characterization of retractively amenable groups that does not refer to nonexpansive actions or retractions?
\end{question}

In view of Proposition~\ref{lem:ra-subgroups}, the author wonders whether $\mathrm{RA}$ can be described, analogously to the original characterization of $\mathrm{EG}$ from \cite{Chou1980}, as the smallest class of groups containing all finite and abelian groups, closed under subgroups, quotients, extensions and direct limits,\footnote{
    As an aside, since all the closure properties in the original characterization of $\mathrm{EG}$ have been established for $\mathrm{RA}$, that characterization could also be used to deduce $\mathrm{EG}\subseteq\mathrm{RA}$.
    Nonetheless, the Chou--Osin characterization provides a more economical argument and, in the proof of Theorem~\ref{thm:bruck-eg-on-retract}, avoids the need to establish subgroup closure for the class $\mathcal M$.
} and satisfying an additional algebraic closure condition that rules out $\mathrm{EG}$ itself as a candidate.

We also ask about the scope of the Bruck extension.

\begin{question}
	\label{q:bruck-beyond-eg}
	Can elementary amenability in Theorem~\ref{thm:bruck-eg-on-retract} be replaced by amenability, or at least by membership in $\mathrm{RA}$?
\end{question}

Additionally, the author considers it worthwhile to investigate elementary amenable analogues of further results for commuting families and invites researchers to pursue this direction.
Examples include the locally convex result in \cite[Theorem~5.3]{DhompongsaKumamSoori2018}, whose extension to amenable semitopological semigroups is explicitly posed as an open problem immediately afterwards.
Further candidates are the multivalued fixed point theorem in \cite[Theorem~3.1]{NananDhompongsa2011} and the fixed point and strong convergence results in \cite[Theorem~3.2, Lemma~4.6 and Theorems~5.1 and~5.4]{Suzuki2006}.

In the follow up, the author plans mostly to study extensions of the presented results to general topologies $\tau$ in the setting of \cite{Borzynski2018}.

\end{document}